\documentclass[12pt]{amsart}
\usepackage{amsmath,amssymb,amscd,amsfonts,verbatim,psfig}

\theoremstyle{definition}

\theoremstyle{remark}

\numberwithin{equation}{section}

\newtheorem{cor}{Corollary}[section]
\newtheorem{prop}[cor]{Proposition}

\theoremstyle{definition}

\newtheorem{defi}[cor]{Definition}

\newtheorem{justnumber}[cor]{}

\newcommand{\co}{\colon\thinspace}    

\newcommand{\grope}{
$$\begin{picture}(140,135)  \footnotesize
    \put(-95,35)       {\psfig{figure=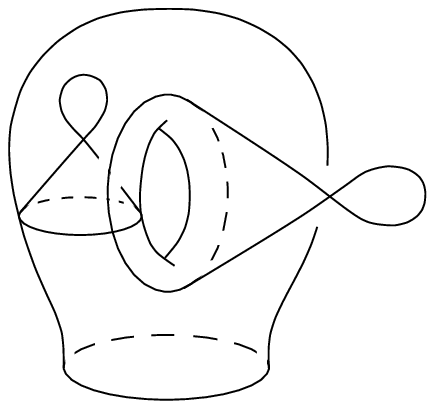}}
    \put(-120,30)       {\psfig{figure=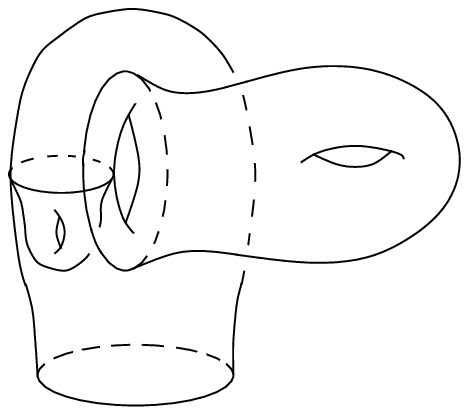}}
\end{picture}$$}

\begin{document}

\title{Surfaces in ${\mathbf 4}$-manifolds and the surgery conjecture}

\author{Vyacheslav S. Krushkal}
\address{Department of Mathematics, University of Virginia, Charlottesville, VA 22904}
\email{krushkal\char 64 virginia.edu}

\thanks{Research was partially supported by NSF grant DMS-0306934}

\date{March 18, 2005}

\begin{abstract}  
We give a survey of geometric approaches to the topological $4$-dimen- sional 
surgery and $5$-dimensional s-cobordism conjectures, with a focus on the study of 
surfaces in $4$-manifolds. The geometric lemma underlying these conjectures
is a statement about smooth immersions of disks and of certain $2$-complexes, 
capped gropes, in a $4$-manifold. 
We also mention a reformulation in terms of the $A-B$ slice 
problem, and the relation of this question to recent developments in the study 
of the classical knot concordance group.
\end{abstract}

\maketitle

\section{Introduction}

The development of the classification theory of topological 
$4$-manifolds in the simply-connected case \cite{F1},\cite{FQ} relied
on the surgery program and parallelled the situation in higher dimensions.
This analogy with the manifolds of dimension greater than $4$
was provided by Freedman's disk embedding theorem \cite{F1} allowing one
to represent hyperbolic pairs in ${\pi}_2 M^4$ by embedded spheres
and therefore to complete surgery. 
Surgery and another ingredient of the classification theory --
the s-cobordism conjecture -- are reduced to closely related 
statements about immersions of disks. (The embedding problem
related to surgery allows certain flexibility in the ambient
manifold, so it is conceivable that the surgery conjecture
holds while the s-cobordism conjecture has an obstruction.)
More recent progress in the subject (cf \cite{FT}, \cite{KQ}), extending 
the class of fundamental groups for which such techniques hold, 
involved the analysis of maps of certain $2$-complexes,
{\em capped gropes}, into $4$-manifolds. Such methods so far have fallen short
of proving the conjectures in general -- the case of free groups
is the key question.
The study of surfaces in $4$-manifolds has been central to the subject, 
and the purpose of this survey is to mention open questions and available 
techniques. In particular, it is worth pointing out that the main open question may
be formulated in terms of the existence of a {\em smooth} immersion, satisfying a 
certain condition on ${\pi}_1$, of a disk 
into a specific smooth $4$-manifold. We also mention another attractive reformulation:
the $A-B$ slice problem which concerns smooth decompositions of the $4$-ball.

{\em Capped gropes}, certain special $2$-complexes embedded in $4$-manifolds,
turned out to be a useful geometric tool for studying homotopies of surfaces
\cite{FQ}, \cite{FT}, \cite{KQ}. One can easily find a capped grope in the 
context of the disk embedding theorem. If one finds a ${\pi}_1$-{\em null} 
caped grope then the problem is solved using the foundational theorem \cite{F1}
that a Casson handle is homeomorphic to the standard $2$-handle. We summarize 
the tools available for manipulating capped gropes in section \ref{definitions}.
On the one hand, these techniques allow one to encode the geometric problem in 
terms of algebra of trees in ${\pi}_1 M$ (and lead to a solution when the 
fundamental group has subexponential growth.) 
On the other hand, this does not immediately lead to an obstruction in the main open case 
(free fundamental group) since the moves on capped gropes, translating the problem into algebra,
capture only a limited class of homotopies of surfaces.

Gropes have recently emerged also as an important tool in a different context: 
the study of the classical knot concordance group \cite{COT}. 
In both applications gropes provide a filtration from bounding
a surface (algebraically: being trivial in homology) to bounding a disk
(being trivial in ${\pi}_1$). However they appear in different guises
(capped vs uncapped) and therefore the techniques for studying them
are quite different. 
Moreover they serve distinct goals: 
in knot theory obstructions measure whether a knot bounds
an embedded grope of a given height. In the surgery context,
one can find a grope of any given height and the problem is
to find a disk. Therefore the question relevant for classification
theory of $4$-manifolds is whether the intersection of the grope filtration 
on links coincides with the class of slice links.

\bigskip 

\centerline{\sc Outline}

1. Formulation of the problem.

2. Symmetric gropes, capped gropes, relation to surgery.

3. The simply-connected case; groups of subexponential growth.

4. Main open problem: canonical examples, free fundamental group.

5. Reformulation of the problem in terms of ${\pi}_1$-null immersions 
$D^2\longrightarrow$ thickening of a capped grope.

6. Reformulation in terms of slicing $Wh(Bor)$ and relation to COT filtration.

7. The $A-B$ slice problem and invariants of decompositions of $D^4$.

\smallskip

{\bf Acknowledgements.} The author is grateful to the organizers of 
the conference on geometry and topology of manifolds at McMaster University 
and of the Banff meeting on knot theory in May 2004. This paper is partially
based on the author's talks at these meetings.

\section{Formulation of the problem} \label{formulation}
 
We start by stating the surgery and s-cobordism 
conjectures. Our main focus will be on the geometry of surfaces
underlying these problems. 

\smallskip

{\bf Surgery conjecture}. {\sl
Let $f\co (M,\partial M)\longrightarrow (X,\partial X)$ be a degree one 
normal map from a topological $4$-manifold to a $4$-dimensional 
Poincar\'{e} duality pair, inducing a homotopy equivalence 
$\partial M\longrightarrow\partial X$. Assume that the surgery obstruction
${\theta}(f)\in L_4({\pi}_1 X)$ vanishes. 
Then $f$ is normally 
bordant to a homotopy equivalence $h\co (N,\partial N)\longrightarrow 
(X,\partial X)$.}

\smallskip

{\bf s-cobordism conjecture}. {\sl
Let $(W^5,M_0,M_1)$ be a $5$-dimensional topological s-cobordism.  
Then $W$ is homeomorphic to the product $W \cong M_0\times I$.}

\smallskip

The proof of both conjectures proceeds as in the higher-dimensional
case, until the following problems remain: in the surgery case 
one has hyperbolic pairs 
$\left( \begin{smallmatrix} 0 & 1\\ 1 & 0 \end{smallmatrix} \right)$
in ${\pi}_2 M$
which need to be represented by embedded spheres. In the s-cobordism
case one can cancel all handles except perhaps for some $2-$ and $3-$ handles,
which are paired up over ${\mathbb Z}{\pi}_1$ by the s-cobordism assumption.
Considering the surgery kernel in the first case, and the cores/cocores of the
handles in the middle level of the cobordism in the second case,
the geometric data is: an immersion of a collection of $S^2\vee S^2$'s into
a $4$-manifold $M$, with a ``distinguished'' intersection point for each pair, 
and with all extra intersection points between the spheres paired up,
with the corresponding Whitney loops contractible in $M$.
In summary, a proof of the following lemma would yield both the surgery
and the s-cobordism conjectures:

\smallskip

{\bf The disk embedding conjecture}. {\sl
Let $(A,{\alpha})\longrightarrow M^4$ be an immersion of a union of disks
with algebraically transverse spheres whose algebraic intersections
and selfintersection numbers are $0$ in ${\mathbb{Z}}[{\pi}_1 M]$.
Then ${\alpha}$ bounds in $M$ disjoint topologically embedded disks with 
the same framed boundary as $A$, and with transverse spheres.}

\smallskip

The {\em algebraically transverse spheres} in the statement 
above are (framed) spheres $\{ S_i\}$ in $M$ with 
$S_i \cdot A_j={\delta}_{i,j}$ (over ${\mathbb{Z}}[{\pi}_1 M]$).
The conjecture is stated for disks; of course to apply it to immersions of a collection
of $S^2\vee S^2$'s above one punctures one sphere in each pair.
Each of the statements above depends on the fundamental group of $M$.
The disk embedding conjecture, proved originally in the simply-connected
case in \cite{F1}, is known to hold for a class of groups including the
groups of subexponential growth \cite{FT},\cite{KQ}. The way that the fundamental
group enters the proof will be clear from the discussion in the following section.
Note that to complete the proof of the surgery conjecture, it would
suffice to solve the disk embedding problem {\em up to s-cobordism}.
There are techniques which use this extra freedom \cite[Chapter 6]{FQ}
but known applications require additional restrictions on the surgery kernel.

\section{Capped gropes and the double point loops.} \label{definitions}

\vspace{.1cm}

In this section we review the 
definitions and terminology for capped gropes and their intersections
in $4$-manifolds. A more detailed exposition may be found in \cite{FQ} or \cite{KQ}.
The section ends with a discussion of the proof of the disk embedding conjecture
in the subexponential growth case and of some related questions.

\begin{figure}[ht]
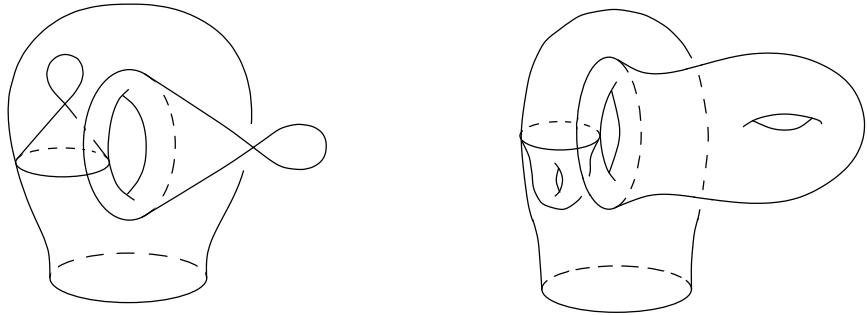
 
\grope
\caption{A capped surface and a grope of height $2$. A {\em capped} grope
of height $2$ is obtained from the $2$-complex on the right by attaching
four caps.}
\label{gropes}
\end{figure}

\begin{defi}
A {\em grope} (or {\em symmetric grope}) $(g,{\gamma})$ is an 
inductively defined pair ($2$-complex, base circle).
A grope has a height $h\in {\mathbb{N}}$. For $h=1$ a grope is a 
compact oriented surface $g$ with a single boundary component $\gamma$. 
A grope $g$ of height $h+1$ is defined inductively: $g$ is obtained
from a grope $g_{h}$ of height $h$ by attaching surfaces (gropes
of height one) to the circles in a symplectic basis for all top stage 
surfaces of $g_{h}$.

A {\em capped grope} $(g^c,{\gamma})$ of height $h$ is obtained from 
a grope $g$ of height $h$ by attaching $2$-cells ({\em caps}) to the 
circles in a symplectic basis for all top stage surfaces of $g$
and introducing finitely many double points among the caps.
Capped gropes of height $1$ are also called capped surfaces, see
figure \ref{gropes}.
\end{defi}

Geometric (and of course algebraic) properties of gropes are substantially different
from those of capped gropes. For a loop $\gamma$ in a space $X$ bounding
a map of a grope of height $n$ is equivalent to being in the $n$-th term
of the derived series of ${\pi}_1 X$. On the other hand, specifying a map of 
a capped grope bounding $\gamma$ means finding special {\em null-homotopies}
of $\gamma$ in $X$. Geometrically, if one considers $X=$ grope of height $n$
then $\gamma$ does not bound a map of a grope of height $n+1$ in $X$.
However given $X=$ capped grope of height $2$, $\gamma$ bounds in $X$ a capped
grope of {\em any} height $\geq 2$. There are several such grope height
raising techniques available \cite{FQ}, \cite{FT}, and the important question
about them is the complexity of the group elements their double point
represent in ${\pi}_1 X$ as a function of height. This is discussed in more detail
below.

Given a capped grope $g^c$, its {\em body} $g$ is the union of all surfaces 
except for the caps. Note that the caps are not allowed to intersect
the body of a capped grope, and the body does not have double
points. This definition is based on the fact
that such $2$-complexes can be found in four-manifolds, bounding Whitney
circles, as they arise in the proof of the surgery and $s$-cobordism
theorems \cite[Theorems 5.1, 7.1, 11.3]{FQ}.

Note that a capped grope $(g^c,{\gamma})$ without double points properly 
embeds into the upper half space $({\mathbb{R}}^3_{+}$, ${\mathbb{R}}^2
\times\{ 0\})$. The ``untwisted'' $4$-dimensional thickening
is obtained from a capped grope $g^c$ without double points by 
first taking its three-dimensional thickening in ${\mathbb{R}}^3_{+}$,
then crossing with $I$, and finally introducing finitely many plumbings
among the caps. The $2$-complex $g^c$ forms a spine of this thickening,
and sometimes we will abuse the notation and use $g^c$ to denote
both objects. To define such ``canonical'' thickening for capped gropes
with double points, one follows the definition above and then introduces
a finite number of plumbings among the thickenings of caps.

Let $(S^c,{\gamma})$ be the untwisted thickening of a capped surface 
(capped grope of height 
$1$) with the body surface $S$ of genus $1$ and with the caps $A$ 
and $B$. There are various ways to get a disk on $\gamma$ in $S^c$, 
we mention some of them to fix the terminology for the rest of the paper. 
Let $({\alpha},{\beta})$ be the symplectic pair of circles in $S$ which 
serve as the attaching curves for $A$ and $B$ respectively. We say that 
$S\smallsetminus ({\alpha}\times I)\hspace{.1cm}\cup\hspace{.1cm}$(two 
parallel copies of $A$) is the {\em surgery on} $S$ {\em along 
the cap} $A$ (analogously one has the surgery on $S$ along $B$.) 
The intermediate operation -- {\em contraction} (or {\em symmetric surgery}) 
-- uses both caps $A$ and $B$, and is described in detail, together with the
associated operation of {\em pushoff}, in \cite[2.3]{FQ}. The definitions 
of surgery and contraction are extended to capped surfaces of higher 
genus, and more generally to capped gropes with more surface stages. 

Let $P$, $Q$ be surfaces or capped gropes in a $4$-manifold $M$. Choose paths 
connecting the base point in $M$ to $P$ and $Q$ to assign groups elements
to the intersection points of the surfaces. Denote by ${\gamma}(p)\in
{\pi}_1 M$ the group element associated to an intersection point 
$p\in P\cap Q$. For self-intersections, one orders the sheets at each 
intersection point to define the group element. Usually the choice of 
an ordering is not specified, since it is either not important, or is 
clear from the context. Given $P$, an immersed surface or a 
capped grope in $M$, to measure the complexity of $P$ in ${\pi}_1 M$
it is important to fix the generators of ${\pi}_1 P$. Assume that 
$P$ is connected. In our applications the source surface (or capped grope) 
will always be simply connected, and we choose the {\em double point 
loops} which pass exactly through one intersection point in the image
(one loop for each self-intersection) as the free generators of ${\pi}_1
P$. Since it is important to 
know how group elements associated to the intersections change under 
the basic operations on surfaces, we briefly review these facts.

\begin{justnumber} \label{surgery and contraction} 
{\bf Surgery, contraction}. 
Suppose $G^c$ is a capped grope, and a surface $P$ intersects one
or both caps in a dual pair of caps $A$, $B$ of $G^c$. Contraction 
of the top stage surface of $G$ to which the caps are attached, or 
surgery along one of the caps doubles the number of intersections, 
but the group elements do not change: if $p'$, $p''$ are the two 
new intersection points created by surgery or contraction from $p$, 
then ${\gamma}(p')={\gamma}(p'')={\gamma}(p)$.
\end{justnumber}

\begin{justnumber} \label{pushoff} {\bf Pushoff}. 
Suppose again that $P$ intersects both caps in a dual pair of caps 
$A$, $B$ of $G^c$, $a\in P\cap A$, $b\in P\cap B$. Pushing $P$ off 
the contraction of $G^c$ along these two caps eliminates the 
intersection points $a$, $b$ and creates two self-intersections
$p'$, $p''$ of $P$. The new group elements are given by
${\gamma}(p')={\gamma}(p'')={\gamma}(a)\cdot{\gamma}(b)^{-1}$.
\end{justnumber}

\begin{justnumber} \label{splitting} {\bf Splitting} \cite{KQ}
is an operation on gropes which arranges all surface stages
above the first one to have genus $1$ and simplifies the pattern of intersections
between the caps at the expense of vastly increasing the genus of the
bottom stage. This operation increases the number of double points
but does not change the group elements represented by them.
\end{justnumber}

\begin{justnumber} \label{dual} {\bf Transverse spheres and gropes}. A {\em 
transverse sphere} 
for a surface $\Sigma$ in $M$ is a framed immersed sphere intersecting 
$\Sigma$ in a single point. Similarly a transverse (capped) grope is a 
(capped) grope whose bottom stage surface intersects $\Sigma$ in a single point.
In a capped grope, the caps and all surfaces, except for the bottom stage,
have a standard transverse grope. More specifically, consider a top stage
surface $S$ of a grope and a pair of dual caps $A,B$, attached along the
curves ${\alpha}, {\beta}$ as above. Consider the normal circle bundle to $S$,
restricted to a parallel copy of $\beta$ (denote it by $T$). 
This is a transverse torus for $A$. Note that a parallel copy of $B$
provides a cap for $T$, and surgering $T$ along this cap gives a 
transverse sphere for $A$. This construction can be continued to produce
transverse gropes of height greater than $1$, cf \cite{KQ}. If a surface $C$ intersects $A$,
taking the sum of $C$ with the capped grope tranverse to $A$ allows one 
to translate between the geometry of intersecting surfaces and the algebra of
the group elements represented by their double point loops (see below). 
\end{justnumber}

The connection to surgery is provided by the following observation.
Under the assumptions of the disk embedding conjecture formulated in 
section \ref{formulation}, the curves $\alpha$ bound
in $M$ capped gropes of height $2$. The proof involves straightforward
manipulations of surfaces \cite[5.1]{FQ}. Given a capped
grope $(G,{\alpha}$) of height two, ${\alpha}$ bounds in it a capped
grope of any given height.
The subtlety is in the complexity of the double point loops
of capped gropes: in the most efficient grope height raising algorithm
presently known \cite{FT} the length of the double point loops of the grope $G_n$ 
of height $n$ constructed in the given grope $G$
grows linearly with $n$ (measured in the free group ${\pi}_1 G$.)

Splitting of a given capped grope (say of height $2$) allows one to 
encode the intersections by trees with edges
labelled by group elements. Vertices of these trees correspond to genus
one pieces of the first stage of the grope, edges correspond to intersections
of caps. For a capped grope of height $2$ this is a $4$-valent tree. 
Splitting uniformizes the intersection patterns of caps -- in particular
all double points of a given cap have the same group element
in ${\pi}_1 M$ associated to them, and this enables one to label the
tree edges by well-defined group elements. In fact the result of 
splitting is stronger, and while locally one sees a tree the global structure 
of this ``intersection graph'' is rather complex.

This algebraic structure 
provides a basis for the proof of the disk embedding theorem for
fundamental groups of subexponential growth in \cite{KQ}: the number of group
elements represented by a tree grows exponentially with height,
while the subexponential growth of ${\pi}_1 M$ allows one to find
trivial group elements. (Recall that a group $G$ generated by $g_1,\ldots, g_k$
has {\em subexponential growth} if, given any ${\alpha}>0$ its growth function 
$gr$ satisfies $gr(n)<e^{{\alpha}n}$ for  large $n$. Here the growth function
$gr\co{\mathbb Z}_{+}\longrightarrow{\mathbb Z}_{+}$ assigns to each positive
integer $n$ the number of elements in $G$ that can be represented as products
of at most $n$ generators $g_{i_1}\cdots g_{i_n}$.)
The translation from algebra back to geometry
of gropes is provided by the moves described above, which yield a capped grope
such that all double point loops are trivial in ${\pi}_1 M$. The proof is
concluded by appealing to the theorem \cite{F1} that a Casson handle is homeomorhic
to the standard $2$-handle.

Such methods have so far fallen short of giving a proof in the general case.
In particular, taking $M^4=$ the untwisted thickening of a capped grope
provides a ``canonical'' problem with free fundamental group:

\smallskip

\begin{justnumber} \label{null}
{\bf ${\mathbf {\pi}_1}$-null disk conjecture}. {\sl Let $(G^c,{\gamma})$ 
denote the untwisted thickening of a $2$-stage capped Grope. 
Then ${\gamma}$ bounds a ${\pi}_1$-null disk in $G^c$.}
\end{justnumber}

\smallskip

Recall that an immersion $f\co D^2\longrightarrow M^4$ is {\em ${\pi}_1$-null}
if the inclusion $f(D^2)\hookrightarrow M$ induces the trivial map on ${\pi}_1$.
This conjecture implies both the surgery and the $s$-cobordism theorems,
since in both contexts one can find a (thickening of a) capped grope.
Note that to complete the proof of the surgery conjecture, one needs
to find a ${\pi}_1$-null disk only up to an s-cobordism of the ambient
$4$-manifold.

We conclude this section with a discussion of the class of ``good'' 
groups for which
the disk embedding conjecture is known to hold. As mentioned above, it is
known to hold for fundamental groups of subexponential growth. Moreover,
it is not difficult to show that the class of groups for which it holds
is closed under extensions and direct limits. It is an interesting question 
whether the disk embedding conjecture holds for amenable groups.
The class of groups containing all groups of subexponential growth and closed
under direct limits and extensions is contained in the class of amenable groups.
The question of whether these two classes actually coincide has been open
until recently: \cite{B} announced a proof that the inclusion is proper and therefore 
gave an example of an amenable group for which the disk embedding conjecture
is open.

\section{Canonical surgery problems, $Wh(Bor)$, and the COT filtration} \label{bor}

Let $Wh(Bor)$ denote the untwisted Whitehead double of the Borromean
rings. This is a $3$-component link -- there is a choice of a $\pm$
clasp for each component, but this choice is not important for the problem.
The slicing question for $Wh(Bor)$ and a related family of links,
with an additional assumption on
 ${\pi}_1$ of the slice complement (see below), 
is a set of canonical surgery problems. (An alternative set 
of ``canonical'' problems is given by the disk embedding conjecture 
in $4$-dimensional thickenings of capped gropes discussed 
above.) Therefore the following question is central to the surgery
program:

\medskip

{\bf Question} \cite{F11} {\sl
Does there exist a $4$-manifold $M$ homotopy equivalent to 
$S^1\vee S^1 \vee S^1$, with the boundary homeomorphic to 
${\mathcal S}^0 (Wh(Bor))$?}

\medskip

The existence of $M$ is equivalent to $Wh(Bor)$ being ``freely slice''
(slice with the additional requirement that the fundamenetal group
of the slice complement in $D^4$ is freely generated by meridians
to the link components). The equivalence is shown by considering
$M=$ slice complement.

Cochran, Orr and Teichner \cite{COT} defined a filtration of 
the classical knot concordance group and introduced an obstruction
theory to show that the quotients of the consecutive terms of
the filtration are non-trivial. A similar result for links
has been established by Harvey \cite{H}. There are two closely
related (and perhaps equivalent) definitions of a filtration,
one of them is in terms of gropes: a knot $K$ is in ${\mathcal F}_n$
if it bounds in $D^4$ an embedded symmetric grope of height $n$.

The following proposition shows that $Wh(Bor)$ lies in the intersection
$\cap {\mathcal F}_n$ of the filtration, and so is not detected
by the aforementioned obstruction theory. Therefore the question,
interesting from the point of view of the surgery conjecture,
is whether the intersection of the filtration coincides with 
the class of slice links.

\smallskip

\begin{prop} 
For any $n$, $Wh(Bor)$ bounds an embedded symmetric
grope of height $n$ in $D^4$.
\end{prop}

\smallskip

\begin{proof}
This proposition holds for a larger class of links (good boundary
links). To be specific we first consider $Wh(Bor)$ and then
we give a more general proof for good boundary links. 

Each component of $Wh(Bor)$
bounds the obvious genus one Seifert surface ${\Sigma}_i, i=1,2,3$ -- 
these are disjoint and lie
in $S^3$. A symplectic basis of curves for these surfaces is given by 
the components $l_i$ of the undoubled link, and a ``short'' curve which bounds
an embedded disk, disjoint from everything else but which is not framed.
(The linking number of this curve with a parallel push-off on the surface
in $S^3$ is $\pm 1$.)
To correct the framing, replace the latter curve by the $(1,\pm 1)$ curve $l'_i$
on the punctured torus. The components of $Bor$
bound disjoint embedded surfaces in $D^4$. More specifically, one of the
components $l_1$ bounds a genus one surface $S$; the other two components
$l_2,l_3$ bound disks $D_2, D_3$. The genus one surface is easily seen in $S^3$; 
we push it in $D^4$ slightly. The two caps of $S$ intersect the disks $D_2$, $D_3$
respectively. Note that $l_i$, $l'_i$ bound disjoint parallel copies of these
respective surfaces ($S,S'$ for $i=1$, $D_i, D'_i$ for $i=2,3$.)

Now assemble the surfaces: $Wh(l_1)$ bounds a $2$-stage capped grope $G_1$
with the base surface ${\Sigma}_1$ (pushed slighly into $D^4$), 
second stage surfaces $S,S'$ and
four caps two of which intersect $D_1$ and $D'_1$, the other two 
intersect $D_2,D'_2$. $Wh(l_2)$, $Wh(l_3)$ bound capped surfaces $G_2,G_3$ with
bases ${\Sigma}_1, {\Sigma}_2$ (pushed into $D^4$),
$D_i, D'_i$ provide {\em embedded} caps
for them. This collection of capped gropes/surfaces doesn't satisfy one 
requirement for being a collection of capped gropes {\em of height} $2$: 
the caps of $G_1$ intersect $D_i, D'_i$ which should be considered as
second stage surfaces for $G_2$, $G_3$. To resolve these interestions, 
note that each second stage disk of ${\Sigma}_2, {\Sigma}_3$ has a
geometrically transverse embedded sphere. To be specific, consider $D_2$.
To see the transverse sphere, consider the transverse capped torus (described 
in \ref{dual}) and surger it along a cap provided by a parallel copy of $D'_2$. 
Now use these transverse spheres to resolve the
intersections between the caps of $G_1$ and the second stage surfaces
of $G_2, G_3$. The result of this construction consists of: embedded disks
bounded by $Wh(l_2)$ and $Wh(l_3)$ and a $2$-stage capped grope bounded
by $Wh(l_1)$. As mentioned in section \ref{definitions}, any of the grope
height raising techniques can be used at this point to find a capped grope
of any height in a neighborhood of $G_1$. Note that the capped gropes we
constructed for $Wh(Bor)$ are of course ${\pi}_1$-null in $D^4$ but 
this is not sufficient to apply the available embedding techniques
(compare with \ref{null}) to conclude that the link is slice. To be useful, the 
nullhomotopies for the double point loops of the caps have to be disjoint from the
body of the gropes -- and finding capped gropes satisfying this requirement
is an open problem.

We will now sketch a proof of the statement for all good boundary links.
By definition here we consider boundary links $L=\partial {\Sigma}$,
${\Sigma}\subset S^3$, such that
there is a basis for $H_1({\Sigma}, {\mathbb Z})$ in which the Seifert pairing 
has the form 
$\left( \begin{smallmatrix} 0 & 1\\ 0 & 0 \end{smallmatrix} \right).$
Examples of good boundary links, important from the perspective of the surgey 
conjecture, are given by the untwisted Whitehead doubles
of links with trivial linking numbers.

Consider a Seifert surface ${\Sigma}$
for $L$ satisfying the good boundary condition, and let ${\alpha}_1,\ldots,
{\alpha}_n,{\beta}_1,\ldots,{\beta}_n$ be curves representing a basis as above.
Consider immersed disks $A_i, B_i$ bounded in $D^4$ by the curves 
${\alpha}_i$, ${\beta}_i$ respectively. Pushing the surfaces ${\Sigma}$ into $D^4$ 
gives a collection of capped surfaces bounded by $L$ in the $4$-ball.
As in \ref{dual} each disk $A_i$ has a transverse torus $T_i$. Consider the sphere $S_i$
given by the surgery on $T_i$ along its cap $B_i$. Using the assumptions
on the Seifert pairing observe that $\{ A_i, S_i\}$ form a collection of disks with 
algebraically transverse spheres, satisfying the assumptions of the disk 
embedding conjecture. Here the $4$-manifold $M$ is the complement of ${\Sigma}$
in $D^4$. Now (as explained in section \ref{definitions}) one finds capped gropes
of height $2$, and therefore of any given height, bounded by $L$ in $D^4$.

\end{proof}

The Borromean rings is the simplest example of a link with trivial linking
numbers but which is homotopically {\em essential} (its components do not bound
disjoint maps of disks in $D^4$.) While the slicing question for $Wh(Bor)$
is a key open problem, it is known \cite{FT} that the untwisted Whitehead doubles
of (a slightly smaller subclass of) homotopically {\em trivial} links are slice.  

It is worth pointing out that while the existence of a $4$-manifold $M$
above ($M\simeq S^1\vee S^1\vee S^1$, ${\partial} M\cong {\mathcal S}^0(Wh(Bor))$)
is a long-standing open problem, one can construct \cite{K1} a double cover $N$ of this
hypothetical manifold $M$ (and of all other canonical manifolds) -- in fact $N$ is
a smooth manifold. Therefore
the surgery conjecture is equivalent to the existence of a free topological 
involution on a certain class of $4$-manifolds. 
Such involutions don't exist {\em smoothly} on at least some of these manifolds,
since by a result of Donaldson \cite{D} the $4$-dimensional surgery conjecture 
fails in the smooth category.

The existence of the manifold $M$ described above and of a closely related family
of $4$-manifolds (equivalently: the free-slice problem for a class of links)
is equivalent to the surgery conjecture. There are no similar canonical problems
presently known specifically for the $s$-cobordism conjecture. A recent paper 
\cite{Q0} develops new techniques in this direction.

\section{Decompositions of $D^4$ and the $A-B$ slice problem}

The $A-B$ slice problem, introduced in \cite{F2}, is a reformulation
of the surgery conjecture. Assume the canonical surgery problem, 
described in section \ref{bor}, has a solution: that is, there exists
a $4$-manifold $M$ homotopy equivalent to $S^1\vee S^1 \vee S^1$, 
with the boundary
homeomorphic to ${\mathcal S}^0 (Wh(Bor))$. It is shown in \cite{F2}
that the compactification of the universal cover $\widetilde M$
is the $4$-ball. The group of covering transformations (the free
group on three generators) acts on $D^4$ with a prescribed
action on the boundary, and roughly speaking the $A-B$ slice 
problem asks whether such action exists -- see \cite{F2} for
a precise formulation. Here we will state an attractive, less
technical formulation, implicitly contained in \cite{F3}.
I would like to thank Michael Freedman for explaining this
approach and for allowing me to include it in this paper.

Let $\mathcal M$ denote the set of (smooth) codimension $0$ 
submanifolds $M$ of $D^4$ with $M\cap \partial
D^4=$ standard $S^1\times D^2\subset \partial D^4$.
We choose a ``distinguished'' curve ${\gamma}\subset \partial M$
whose neighborhood is the solid torus specified above.

Consider invariants of such submanifolds:
$I:{\mathcal M}\longrightarrow \{0,1\}$ satisfying axioms 1-3:

\smallskip

{\sl Axiom 1.} $I$ is a topological invariant: if $(M, {\gamma})$ is
diffeomorphic to $(M', {\gamma}')$ then 
$I(M,{\gamma})=I(M',{\gamma}')$.

\smallskip

{\sl Axiom 2.} If $(M,{\gamma})\subset (M',{\gamma}')$ and $I(M,{\gamma})=1$
then $I(M',{\gamma}')=1$.

\smallskip

For a codimension $0$ submanifold $A$ of $D^4$ let $\partial A={\partial}^{-}A
\cup{\partial}^{+}A$ where ${\partial}^{+}A=\partial A
\cap \partial D^4$.
We say that $D^4=A\cup B$ is a {\em decomposition} of $D^4$ 
if $A,B$ are codimension zero submanifolds, and $\partial D^4=
\partial^{+}A\cup \partial^{+}B$ is the standard genus one Heegaard
decomposition of $S^3$ (of course $\partial^{-}A= \partial^{-}B$.)

\smallskip

{\sl Axiom 3.} If $(A,B)$ is a decomposition of $D^4$ then $I(A)+I(B)=0$.

\smallskip

Note that $I$ is an invariant of a submanifold and does not depend on an
embedding $M\hookrightarrow D^4$. One can require an invariant to be defined
on the class of all pairs ($4$-manifold $M$, distinguished circle in $\partial M$)
but only submanifolds of $D^4$ are relevant for the surgery conjecture.

There is an elementary example of such an invariant (in any dimension)
given by homology: 
define $I_h(M,{\gamma})=0$ or $1$ depending on whether ${\gamma}\neq 0$
or ${\gamma}=0$ respectively in $H_1(M)$ (with any fixed 
coefficients). Axioms 1 and 2 are satisfied automatically and axiom
3 follows from Alexander duality.

\medskip

Are there any other invariants? In particular, the property connecting
this question to $4$-dimensional topology and which is relevant 
for the $(A,B)$-slice problem is summarized in the additional axiom stated
below.
Given $(M',{\gamma}'), (M'',{\gamma}'') \in{\mathcal M}$, define
the ``double'' $D(M',M'')=(S^1\times D^2\times I)\cup (M'\cup M'')$ 
where $M',M''$ are attached to $S^1\times D^2\times \{1\}$ along
the distinguished solid tori in their boundaries, so that 
${\gamma}',{\gamma}''$ form the Bing double of the core of the solid torus.
(Note that if both $M'$, $M''$ embed in $D^4$ then so does $D(M',M'')$,
and therefore $D(M', M'')$ is an element of ${\mathcal M}$.)
Let $\gamma$ denote the core of $S^1\times D^2\times\{ 0\}$.

\smallskip

{\sl Axiom 4.} Let $(M',{\gamma}'), (M'',{\gamma}'') \in{\mathcal M}$ be
such that $I(M',{\gamma}')=I(M'',{\gamma}'')=1$. Then
$I(D(M',M''),{\gamma})=1$.

\smallskip

The homology candidate above clearly doesn't satisfy this last axiom:
consider $(M',{\gamma}')=(M'',{\gamma}'')=(D^2\times D^2, \{0\}\times 
\partial D^2)$. Then $I_h=1$ for both $M', M''$. However $D(M',M'')$
is obtained from the collar $S^1\times D^2\times I$ by attaching two
$2$-handles along the Bing double of the core of the solid torus,
so this core is not trivial in homology of $D(M',M'')$ and $I_h(D(M',M''))=0$.

An invariant of decompositions satisfying axioms 1-4 would be a very good
candidate for an obstruction to surgery: the connection is made by considering
fundamental domains of a hypothetical action of the free group on $D^4$,
discussed above.
A new approach to constructing an invariant of decompositions is outlined
in \cite{K2}. That paper defines an invariant of $4$-manifolds using
the notion of link-homotopy; in a certain sense it is designed to 
satisfy axiom 4. It would provide an obstruction to surgery if it satisfied
``Alexander duality'' (axiom 3).

\vspace{.2cm}

\end{document}